\documentclass{article}
\usepackage{amssymb}

\let\le=\leqslant
\let\ge=\geqslant

\def\mup#1{{(-1)}^{#1}}
\def\mun{\mup n}
\def\ap#1{{(\alpha\,{+}\,#1)}}\def\bp#1{{(\beta \,{+}\,#1)}}
\def\apn{\ap n}
\def\Li{\mathop{\mathrm{Li}}\nolimits}
\def\Re{\mathop{\mathrm{Re}}\nolimits}

\def\bigstrut{{\strut\over\strut}}
\def\Bigstrut{{\bigstrut\over\strut}}

\def\Cset{\mathbb{C}}

\begin{document}

\title{On a multiple harmonic power series.}

\author{Michel \'Emery
\footnote{IRMA, 7 rue Ren\'e Descartes, 
      67084 Strasbourg Cedex, France.
     \rm\texttt{emery@math.u-strasbg.fr}}}
     
\maketitle

\begin{abstract}
If $\Li_s$ denotes the polylogarithm of order~$s$, 
where $s$ is a natural number, and if $z$ belongs
to the unit disk,
\[
  \Li_s {\Bigl({-z\over1\,{-}\,z}\Bigr)}
  =-\sum _{1\le i_1\le \ldots\le i_s}
  {z^{i_s}\over i_1i_2\ldots i_s}\;.
\]
In particular,
\[
  \sum _{n\ge 1}{\mup {n+1}\over n^s\!}\,=
  \sum _{1\le i_1\le \ldots\le i_s}
  {1\over{i_1}\ldots{i_s}\,2^{i_s}}\;.
\]
\bigbreak
\noindent
\textbf{MSC}: 11M35, 33E20, 40A25, 40B05.

\medbreak
\noindent
\textbf{Keywords}: Multiple harmonic series, 
Lerch function, Polylogarithm.

\end{abstract}

\bigbreak
\goodbreak\noindent\textbf{Introduction}\\

Equalities and identities between multiple 
harmonic series and polylogarithms 
have been investigated by many authors; see for 
instance~\cite{DB:PB:DJB:PL}
and the references therein. 
These series usually involve summations 
over all $s$-tuples ${(i_1,\ldots,i_s)}$ of 
natural numbers such that ${i_1<\ldots<i_s}$,
where $s$~is fixed. We shall be concerned with an
instance  where the summation indices may be 
equal to each other, that is, a sum over all 
integers ${i_1,\ldots ,i_s}$ verifying 
${1\le i_1\le \ldots\le i_s}$.\\

\goodbreak\noindent\textbf{Definition.} 
\it For $\,\alpha \in \Cset
\setminus{\{-1,-2,\ldots\}}$ and $\,s\in {\{1,2,\ldots\}}$, 
the Lerch function of order~$\,s$ with 
shift~$\,\alpha $ is defined by\rm
\[
\Li_s^\alpha (w)=\sum _{n\ge 1}{w^n\over\apn^s\!}\,\;.
\]

The power series converges in the unit disk only,
but the analytic function $\Li_s^\alpha $ 
extends to the whole complex plane
minus a cut along for instance the half-line ${[1,\infty )}$;
see \S~1.11 of~\cite{E:HTF}. 
When ${\alpha=0}$, $\Li^\alpha_s$ is 
just $\Li_s$, the usual
polylogarithm of order~$s$.
We shall be interested in the values of 
$\Li_s^\alpha (w)$ in the half-plane 
${\Re (w)<{1\over2}}$ only; 
they are given by the next proposition (which 
establishes anew the existence of the
analytic extension of $\Li^\alpha_s$
to that half-plane).\\

\goodbreak\noindent\textbf{Proposition.} 
\it For $\,\alpha \in 
\Cset\setminus{\{-1,-2,\ldots\}}$ and
$\,s\in {\{1,2,\ldots\}}$, one has
the following power series expansion,
which converges for $\,{|z|<1}$:\rm
\begin{eqnarray}
\Li_s^\alpha {\Bigl({-z\over1\,{-}\,z}\Bigr)}&\nonumber\\
&\displaystyle\hskip-13mm 
=-\hskip-3mm \sum _{1\le i_1\le \ldots\le i_s}
{(i_s\,{-}\,1)!\over
(\alpha {+}1)(\alpha +2)\ldots(\alpha +i_s)}\ \>
{z^{i_s}\over
(\alpha {+}i_1)(\alpha +i_2)\ldots(\alpha +i_{s-1})}\;.
\nonumber\end{eqnarray}
\\

When ${\alpha =0}$, the right-hand side becomes much 
simpler:\\

\goodbreak\noindent\textbf{Corollary.} 
\it Fix $\,s\in {\{1,2,\ldots\}}$. For
$\,z$ in the unit disk,
\[
\Li_s {\Bigl({-z\over1\,{-}\,z}\Bigr)}
=-\sum _{1\le i_1\le \ldots\le i_s}
{z^{i_s}\over i_1i_2\ldots i_s}\;.
\]

Equivalently, for $\,{\Re w<{1\over2}}$,\rm
\[
\Li_s {(w)}
=-\sum _{1\le i_1\le \ldots\le i_s}
{1\over i_1i_2\ldots i_s}\>
{\Bigl({-w\over1\,{-}\,w}\Bigr)}^{i_s}\;.
\]

For instance, choosing ${z={1\over2}}$ in the corollary gives
\[
\sum _{n\ge 1}{\mup {n+1}\over n^s\!}\,=
\sum _{1\le i_1\le \ldots\le i_s}
{1\over{i_1}\ldots{i_s}\,2^{i_s}}\;;
\]
as the left-hand side is ${(1\,{-}\,2^{1-s})\zeta (s)}$, 
and the right-hand one is ${\sum a_p/(p2^p)}$ with 
${a_p\le {(1\,{+}\,\log p)}^{s-1}}$, this is a reasonably 
fast series expansion of~$\zeta (s)$. More generally, 
taking ${z={1\over2}}$ in the proposition gives a
more rapidly convergent series for 
${\sum\limits_n{(-1)}^n{(\alpha\,{+}\,n)}^{-s}}$.\\

\goodbreak\noindent\textbf{Proof of the proposition}\\

Observe first that the right-hand side in the statement of
the proposition is a power series 
$\sum c_pz^p$, with
\begin{eqnarray}
{|c_p|}&\displaystyle\le {(p\,{-}\,1)!
\over{|\alpha {+}1|}\ldots{|\alpha {+}p|}}
\sum _{1\le i_1\le \ldots\le i_{s-1}\le p}
{1\over{|\alpha {+}i_1|}\ldots
{|\alpha {+}i_{s-1}|}}\nonumber\\
&\displaystyle\hskip1cm 
\le {(p\,{-}\,1)!\over{|\alpha {+}1|}\ldots{|\alpha {+}p|}}
\ {p^{s-1}\over C{(\alpha )}^{s-1}}\;,
\nonumber\end{eqnarray}
where $C(\alpha )=\inf{\{{|\alpha {+}1|},
{|\alpha +2|},\ldots\}}>0$.
By d'Alembert's test, the series $\sum {|c_pz^p|}$
converges in the unit disk; so does also $\sum c_pz^p$.
Consequently, to prove the proposition, it suffices 
to show that both sides of the claimed identity 
are equal for ${|z|}$ small enough. We shall take 
${|z|<{1\over2}}$; for such a~$z$,
one has ${|{-z/{(1{-}z)}}|<1}$, whence
\[
\Li_s^\alpha {\Bigl({-z\over1\,{-}\,z}\Bigr)}
=\sum _{n\ge 1}{\mun\over\apn^s}
{\Bigl({z\over1\,{-}\,z}\Bigr)}^n
=\sum _{n\ge 1}{\mun\over\apn^s}
\sum _{p\ge n}{\Bigl({p{-}1\atop n{-}1}\Bigr)}z^p\;.
\]
Exchanging the summations gives
$$
\Li_s^\alpha {\Bigl({-z\over1\,{-}\,z}\Bigr)}
=\sum _{p\ge 1}z^p\sum _{n=1}^p
{\Bigl({p{-}1\atop n{-}1}\Bigr)}{\mun\over\apn^s}\;.
\leqno(*)
$$
This exchange is licit because, setting ${K(\alpha )=1/
\smash{\inf\limits_{n\ge 1}{|1{+}\alpha /n|}}<\infty }$,
one has the estimate
\[
\!\sum _{p\ge n\ge 1}{\biggl|{\Bigl({p{-}1\atop n{-}1}\Bigr)}
{\!\mun {z}^p\over{(\alpha {+}n)}^s}\biggr|}
\le K{(\alpha )}^s
\!\!\sum _{p\ge n\ge 1}\!{\Bigl({p{-}1\atop n{-}1}\Bigr)}
{{|z|}^p\over n^s}=K{(\alpha )}^s\!\sum _{n\ge 1}{1\over n^s}
{\Bigl({{|z|}\over1{-}{|z|}}\Bigr)}^{\!n},
\]
which is finite because ${0\le {|z|}/(1\,{-}\,{|z|})<1}$.

To establish the proposition, it remains to compute the 
coefficient of~$z^p$ in~$(*)$, and more precisely to show 
that, for all ${\alpha \notin {\{-1,-2,\ldots\}}}$,
${s\in {\{1,2,\ldots\}}}$ and ${p\ge 1}$,
\[
\sum _{n=1}^p
{\Bigl({p{-}1\atop n{-}1}\Bigr)}{\mun\over\apn^s}=
\left\{\begin{array}{c}
   \vphantom{\Bigstrut}\end{array}\right.\mkern-25mu
   \begin{array}{ll}
     \displaystyle
     -{(p\,{-}\,1)!\over\ap1_p}\;
     \sum _{1\le i_1\le \ldots\le i_{s-1}\le p}
     \;\prod _{r=1}^{s-1}{1\over\alpha \,{+}\,i_r}
     &\mbox{ if $s\ge 2$,}\\
     \noalign{\vskip2pt} \displaystyle
     -{(p\,{-}\,1)!\over\ap1_p}
     &\mbox{ if $s=1$,}\\
  \end{array}
\]
where ${(x)}_p$ is the Pochhammer symbol 
standing for
${\prod\limits _{j=0}^{p-1}(x\,{+}\,j)}$.

Putting ${q=p\,{-}\,1}$, ${m=n\,{-}\,1}$ and
${\beta =\alpha \,{+}\,1}$, 
it suffices to prove the following lemma:\\

\goodbreak\noindent\textbf{Lemma.} 
\it For all $\,{\beta \in \Cset\setminus
{\{0,-1,\ldots\}}}$
and all integers $\,{q\ge 0}$ and $\,{s\ge 1}$, one~has
\[
\sum _{m=0}^q
{\Bigl({q\atop m}\Bigr)}{\mup m\over\bp m^s}=
\left\{\begin{array}{c}\vphantom{{\strut\over\strut}
\over\strut}\end{array}\right.\mkern-25mu
   \begin{array}{ll}
     \displaystyle
     {q!\over{(\beta )}_{q+1}}\;
     \sum _{0\le i_1\le \ldots\le i_{s-1}\le q}
     \;\prod _{r=1}^{s-1}{1\over\beta \,{+}\,i_r}
     &\mbox{ if $s\ge 2$;}\nonumber\\
     \noalign{\vskip2pt}\displaystyle
     {q!\over{(\beta )}_{q+1}}
     &\mbox{ if $s=1$.}\nonumber
   \end{array}
\]\rm

\medskip

\goodbreak\noindent\textbf{Proof of the lemma}\\

Call $L(q,\beta )$ and $R(q,\beta )$ 
the left- and right-hand
sides in this statement.

First, for ${q=0}$ and ${q=1}$, the lemma is 
easily verified:
\begin{eqnarray}
&\displaystyle L(0,\beta )=
{1\over{\mkern0mu\beta}^s}=R(0,\beta )\;;
&\nonumber\\
\noalign{\vskip2pt}
&\displaystyle L(1,\beta )=
{1\over{\mkern0mu\beta}^s}-{1\over\bp1^s}
={1\over\beta \bp1}
\sum _{u+v=s-1}{1\over{\mkern0mu\beta}^u}\>
{1\over\bp1^v}=R(1,\beta )\;.&
\nonumber\end{eqnarray}
Next, observe that
\begin{eqnarray}
L(q{+}1,\beta )&=&
\sum _m{\Bigl({q{+}1\atop m}\Bigr)}
{\mup m\over\bp m^s}
=\sum _m{\Bigl[
{\Bigl({q\atop m}\Bigr)}+{\Bigl({q\atop m{-}1}\Bigr)}
\Bigr]}{\mup m\over\bp m^s}\nonumber\\
&=&\sum _m
{\Bigl({q\atop m}\Bigr)}{\Bigl[
{\mup m\over\bp m^s}+{\mup {m+1}\over\bp {m\,{+}\,1}^s}
\Bigr]}\nonumber\\
\noalign{\vskip5pt}&=&L(q,\beta )-L(q,\beta {+}1)\;.
\nonumber\end{eqnarray}
The rest of the proof will consist in verifying that,
for ${q\ge 1}$, 
the right-hand side satisfies the same relation:
$$
R(q{+}1,\beta ) = R(q,\beta )-R(q,\beta {+}1)\;.\leqno(**)
$$
When this is done, $(**)$ immediately implies 
that the property
\[
\forall \beta \notin {\{0,-1,\ldots\}}\qquad
L(q,\beta )=R(q,\beta )
\]
extends by induction from ${q=1}$ to all ${q\ge 2}$, 
thus proving the lemma.

The following notation will 
simplify the proof of~$(**)$:
for all integers ${n\ge 0}$, ${t\ge 1}$, 
${a\ge 0}$ and ${b\ge a}$, set
\[
f_n={1\over\beta \,{+}\,n}\;;\qquad S_a^b(t)=
\sum _{a\le i_1\le \ldots\le i_{t}\le b}
\;\prod _{r=1}^{t}f_{i_r}\;;\qquad S_a^b(0)=1\;;
\]
and remark that, for 
${t\ge 0}$ and ${0\le a\le b<c}$,
\[
S_a^c(t)=\sum _{u+v=t}S_a^b(u)\,S_{b+1}^c(v)\;;
\]
similarly, for ${1\le a\le b}$,
\[
S_{a-1}^{b+1}(t)=\sum _{u+v+w=t}
S_{a-1}^{a-1}(u)\,S_a^b(v)\,S_{b+1}^{b+1}(w)
=\sum _{u+v+w=t}
f_{a-1}^u\,S_a^b(v)\,f_{b+1}^w\;.
\]
Keeping these remarks in mind, the proof of $(**)$
goes as follows. Put ${t=s\,{-}\,1}$ and write
\begin{eqnarray}
R(q,\beta )-R(q,\beta {+}1)&=&
{q!\over{(\beta )}_{q+1}\!\!}\,\,S_0^q(t)
-{q!\over{(\beta {+}1)}_{q+1}\!\!}\,\,S_1^{q+1}(t)
\nonumber\\
&=&{q!\over{(\beta {+}1)}_q\!\!}\,\,
{\Bigl[{1\over\beta }\,S_0^q(t)-
{1\over\beta {+}q{+}1}\,S_1^{q+1}(t)\Bigr]}\;.
\nonumber\end{eqnarray}
The quantity in square brackets can be rewritten as 
\begin{eqnarray}
\lefteqn{f_0\sum _{v=0}^{t}
S_0^0(t{-}v)S_1^q(v)
-f_{q+1}\sum _{v=0}^{t}S_1^q(v)
S_{q+1}^{q+1}(t{-}v)}\nonumber\\
&\qquad\quad=&
\sum _{v=0}^{t}S_1^q(v){(f_0^{t-v+1}-f_{q+1}^{t-v+1})}
=\sum _{v=0}^{t}S_1^q(v)\,{(f_0\,{-}\,f_{q+1})}\!\!
\sum _{u+w=t-v}\!\!f_0^uf_{q+1}^w\nonumber\\
&\qquad\quad=&
{(f_0\,{-}\,f_{q+1})}\!\!\sum _{u+v+w=t}\!\!
f_0^uS_1^q(v)f_{q+1}^w
={(f_0\,{-}\,f_{q+1})}\,S_0^{q+1}(t)\nonumber\\
&\qquad\quad=&
{q+1\over\beta (\beta {+}q{+}1)}S_0^{q+1}(t)\;.
\nonumber\end{eqnarray}
Finally, the difference 
${R(q,\beta )-R(q,\beta {+}1)}$ 
amounts to
\[
{q!\over{(\beta {+}1)}_q\!\!}\ \>
{q+1\over\beta (\beta {+}q{+}1)}\,S_0^{q+1}(t)
={(q\,{+}\,1)!\over{(\beta )}_{q+2}\!\!}\,\,
S_0^{q+1}(t)=R(q{+}1,\beta )\;.
\]
This proves $(**)$, and at the same time the lemma 
and the proposition.\\

\goodbreak\noindent\textbf{Remark.} In the
particular case when ${\alpha =0}$ and 
${z={1\over2}}$, $(*)$ is exactly formula~(4)
of J.~Sondow~\cite{JS}, obtained by
accelerating convergence of the 
alternating zeta series. What the lemma
does is computing the numerator in 
that formula. To that end, $(**)$ is needed
for ${\beta \in \{1,2,\ldots\}}$ only;
but the general proof is just as easy.

\end{document}